\documentclass[nrl,nonblindrev]{informs3}
\usepackage{optidef}
\usepackage{IEEEtrantools}
\usepackage{booktabs}
\usepackage{multirow}
\usepackage{amsmath}
\usepackage{amsfonts}
\usepackage{caption}
\usepackage{graphicx}
\usepackage{subcaption}
\usepackage{algorithm, algpseudocode}
\usepackage{pgfplots}
\pgfplotsset{compat=1.8}
\usepgfplotslibrary{statistics}
\usepackage{pgf} 
\usepackage{pgfplots}
\usetikzlibrary{fit,calc}
\usepgfplotslibrary{external}
\usetikzlibrary{decorations.pathmorphing}
\usetikzlibrary{patterns}

\usepackage{rotating}
\DoubleSpacedXI % Made default 4/4/2014 at request
\usepackage{endnotes}
\let\footnote=\endnote

\usepackage{natbib}
 \bibpunct[, ]{(}{)}{,}{a}{}{,}%
 \def\bibfont{\small}%
 \def\BIBand{and}%

\TheoremsNumberedThrough     % Preferred (Theorem 1, Lemma 1, Theorem 2)
\ECRepeatTheorems

\EquationsNumberedThrough    % Default: (1), (2), ...
\MANUSCRIPTNO{xxxx} % When the article is logged in and DOI assigned to it,
\begin{document}
%%%%%%%%%%%%%%%%

% Outcomment only when entries are known. Otherwise leave as is and
%   default values will be used.
%\setcounter{page}{1}
%\VOLUME{00}%
%\NO{0}%
%\MONTH{Xxxxx}% (month or a similar seasonal id)
%\YEAR{0000}% e.g., 2005
%\FIRSTPAGE{000}%
%\LASTPAGE{000}%
%\SHORTYEAR{00}% shortened year (two-digit)
%\ISSUE{0000} %
%\LONGFIRSTPAGE{0001} %
%\DOI{10.1287/xxxx.0000.0000}%

% Author's names for the running heads
% Sample depending on the number of authors;
% \RUNAUTHOR{Jones}
% \RUNAUTHOR{Jones and Wilson}
% \RUNAUTHOR{Jones, Miller, and Wilson}
% \RUNAUTHOR{Jones et al.} % for four or more authors
% Enter authors following the given pattern:
\RUNAUTHOR{Kiani, Isik, Eksioglu}

% Title or shortened title suitable for running heads. Sample:
% \RUNTITLE{Bundling Information Goods of Decreasing Value}
% Enter the (shortened) title:
\RUNTITLE{Dynamic Tuberculosis Screening for Healthcare Employees}

% Full title. Sample:
% \TITLE{Bundling Information Goods of Decreasing Value}
% Enter the full title:
\TITLE{Dynamic Tuberculosis Screening for Healthcare Employees}

% Block of authors and their affiliations starts here:
% NOTE: Authors with same affiliation, if the order of authors allows,
%   should be entered in ONE field, separated by a comma.
%   \EMAIL field can be repeated if more than one author
\ARTICLEAUTHORS{%
\AUTHOR{Mahsa Kiani$^{a}$ , Tugce Isik$^{a}$, Burak Eksioglu$^{b}$}
\AFF{$^{a}$Department of Industrial Engineering, Clemson University, Clemson, SC 29634 \\ $^{b}$Department of Industrial Engineering, University of Arkansas, Fayetteville, AR 72701}
\AFF{\EMAIL{mkiani@clemson.edu}, \EMAIL{tisik@clemson.edu}, \EMAIL{burak@uark.edu}}

% Enter all authors
} % end of the block

\ABSTRACT{Regular tuberculosis (TB) screening is required for healthcare employees since they can come into contact with infected patients. TB is a serious, contagious, and potentially deadly disease. Early detection of the disease, even when it is in latent form, prevents the spread of the disease and helps with treatment. Currently, there are two types of TB diagnostic tests on the market: skin test and blood test. The cost of the blood test is much higher than the skin test. However, the possibility of getting a false positive or false negative result in skin test is higher especially for persons with specific characteristics, which can increase costs. In this study, we categorize healthcare employees into multiple risk groups based on the department they work in, the specific job they do, and their birth country. We create a Markov decision process (MDP) model to decide which TB test should be taken by each employee group to minimize the total costs related to testing, undetected infections, employees' time lost.  Due to the curse of dimensionality, we use approximate dynamic programming (ADP) to obtain a near-optimal solution. By analyzing this solution to the ADP we specify not only the type but the frequency with which each test should be taken. Based on this analysis, we propose a simple policy that can be used by healthcare facilities since such facilities may not have the expertise or the resources to develop and solve sophisticated optimization models.}

% Sample
%\KEYWORDS{deterministic inventory theory; infinite linear programming duality;
%  existence of optimal policies; semi-Markov decision process; cyclic schedule}

% Fill in data. If unknown, outcomment the field
\KEYWORDS{Markov decision process, dynamic programming, tuberculosis infection, tuberculosis testing, healthcare employee screening}

\maketitle
%%%%%%%%%%%%%%%%%%%%%%%%%%%%%%%%%%%%%%%%%%%%%%%%%%%%%%%%%%%%%%%%%%%%%%

% Samples of sectioning (and labeling) in OPRE
% NOTE: (1) \section and \subsection do NOT end with a period
%       (2) \subsubsection and lower need end punctuation
%       (3) capitalization is as shown (title style).
%
%\section{Introduction.}\label{intro} %%1.
%\subsection{Duality and the Classical EOQ Problem.}\label{class-EOQ} %% 1.1.
%\subsection{Outline.}\label{outline1} %% 1.2.
%\subsubsection{Cyclic Schedules for the General Deterministic SMDP.}
%  \label{cyclic-schedules} %% 1.2.1
%\section{Problem Description.}\label{problemdescription} %% 2.
% Text of your paper here

   \section{Introduction}
Tuberculosis (TB) is an infectious disease caused by Mycobacterium tuberculosis (Mtb) that mainly affect the lungs, but can also impact other parts of the body. Approximately one third of world's population is infected by Mtb \cite{getahun2015latent}. TB infections are categorized into two groups: latent and active. Individuals with latent infection have the disease but do not show any symptoms, and they are not infectious. However, a latent infection can turn in to active TB. Individuals with active TB are infectious and spread TB. Thus, medication is required for active TB. Although TB infections have been declining in Western Europe and North America, there are groups of people with high prevalence of TB such as immigrants and prisoners \cite{youakim2016occupational}. Due to being in contact with these groups, healthcare providers face a higher than usual risk of exposure to TB. Early detection of TB infections is critical to control the spread of the disease. Thus, healthcare workers are usually suggested to take a TB diagnosis test upon initial hiring and thereafter at regular intervals. In fact, CDC guidelines prior to 2005 suggested annual testing for medium or higher risk environments. However, these guidelines were later changed to allow for other policies to be implemented. Yet, many healthcare facilities such as Prisma Health still follow the old guidelines which usually result in expensive testing for such facilities. 

There are currently two common tests for detection of latent and active TB on the market: the skin test and the blood test which is also called Interferon Gamma Release Assay (IGRA). The skin test involves an injection of 0.1 mL of a liquid containing five tuberculin units of purified protein derivative (PPD) into the top layers of skin of the forearm. Once this liquid is injected, the test has to be read within 48-72 hours. A confirmed positive test involves swelling of the injection site, while a negative test has no signs of inflammation. Although the skin test is rapid, it can be prone to false positive results that require expensive further testing. There are different factors that may cause a false positive result such as sensitization of the test to some nontuberculous mycobacteriaes, incorrect interpretation of reaction, and incorrect method of the skin test administration \cite{menzies2007new,youakim2016occupational}. The skin test also shows a false positive if the patient has already had the Bacillus Calmette-Gurin (BCG) vaccine. False positive results cause  administration of unnecessary chest X-rays at an extra cost of \$100-\$400 per X-ray \cite{Chest:2020}. There are also factors such as incorrect interpretation of reaction, some viral illnesses (e.g., measles and chicken pox) and recent live-virus vaccination (e.g., measles and smallpox) that may lead to a false negative skin test result \cite{youakim2016occupational}. Individuals who received a false negative skin test may infect others until they begin showing obvious symptoms and begin treatment. 

In a TB blood test, a small amount of blood is drawn and sent to a laboratory. Thus, healthcare employees visit the clinic only one time, as opposed to twice for a skin test. The cost of a blood test is higher than that of a skin test, but the test is more accurate. Also, unlike a skin test, the accuracy of a blood test is not affected by prior BCG vaccinations. Although the cost of performing a blood test is much higher than the cost of a skin test, blood tests can still be preferable at least for specific groups of employees when considering the lost time of employees and healthcare professionals who administer the test, as well as the costs related to false positive or false negative results.

In this study, we collaborate with Prisma Health, a healthcare system in South Carolina, that currently requires two skin tests for new employees and an annual skin test for all other employees. Prisma Health does not utilize blood tests currently. As part of our analysis, we categorize the employees into multiple groups based on the risk of infection related to their job, their work environment, birth country and BCG vaccination history. We define an infection rate for each employee group which depends on the number of infected people the employees in that group potentially get in contact with. We introduce a Markov Decision Process (MDP) model to develop a screening plan for a healthcare facility by determining the type and frequency of TB test to be used for each employee group. The objective is to minimize the expected total cost of the system. Due to the curse of dimensionality, we use Approximate Dynamic Programming (ADP) to find a ``near-optimal" solution. Based on this solution, we propose a benchmark policy that is easy to implement by the healthcare facility, and evaluate this simple policy using data obtained from Prisma Health.

\section{Literature Review} 
There are two streams of literature related to our research. We briefly explain each and highlight our contributions in relation to other studies.

The first stream of studies relates to the analysis of TB screening and its importance for healthcare employees. As mentioned in the previous section, healthcare employees are at a higher than usual risk of getting TB infections. \cite{hazard2016hidden} show that even healthcare employees who work in places that are not in direct contact with patients, for instance employees who work in a hospital kitchen, are at a higher risk of being infected. Most healthcare facilities have a regular plan of TB screening for their employees. One commmon strategy is to administer an annual skin test. \cite{mullie2017revisiting} tested three different scenarios to propose a TB screening program. The first scenario proposes annual screening for all employees. In the second scenario, only employees with high-risk tasks, such as respiratory therapy, are tested yearly and other employees are tested only after recognized exposure. The third scenario tests all employees only after recognized exposure. They evaluated these scenarios by using both skin and blood tests. They also did a cost effectiveness study on 1000 US healthcare employees with no positive TB history. Results of their experiments indicated that for most US healthcare employees annual TB testing is expensive with limited health gains. Thus, regular annual testing may not be an effective strategy for most health systems.

In addition to the testing frequency, the type of TB diagnosis test to be used is also an important decision. \cite{talbot2012specificity} studied the specificity of skin and blood tests among students in a low-tuberculosis incidence setting. They concluded that the blood test performed better in these settings.  \cite{khalil2013comparison} performed a sensitivity comparison of blood and skin tests in 50 cases with active TB and showed that the sensitivity of the blood test is about 80$\%$ while the skin test sensitivity (accuracy ) is about 28$\%$.  Accuracy of blood and skin tests in detecting latent and active TB is  also presented in other studies \cite{leung2010t, foster2014tuberculosis, leung2015t, king2015t}. However, further cost analysis studies are needed in this area.
\cite{de2009cost} performed a cost-effectiveness analysis of TB blood and skin tests by considering the direct test costs and cost of missed work time. \cite{wrighton2012screening} extended this study and considered the performance of each test in calculation of total cost of taking a blood test versus a skin test. Thus, cost of subsequent tests and treatments in the case of getting a positive result was included in the total cost. Neither of these studies considered the potential costs that a false negative result may cause. Previous studies  also did not analyze the cost-effectiveness of the tests with respect to employee characteristics. 

The second stream of relevant studies is on the application of MDP models in the context of prevention, screening, and treatment of diseases. These decisions are typically made sequentially over long periods in uncertain environments \cite{ni2017markov}. In addition to the patient's current health status, the uncertainty in progression of the disease, impact of the treatment on the patient, and accuracy of the test results have to be considered in determining the treatment decision \cite{steimle2017markov}. Using MDPs is often appropriate to analyze such problems since the decisions are made sequentially over time in a fundamentally stochastic environment. \cite{eghbali2019markov} developed a MDP to model adverse drug reactions in medication treatment of type 2 diabetes. MDPs are also used in breast cancer screening \cite{maillart2008assessing, chhatwal2010optimal, ayer2012or, cevik2018analysis}, treatment of HIV \cite{shechter2008optimal}, and public policy decisions related to the transmission of communicable diseases \cite{kaplan1995probability, zaric2001optimal}. \cite{ayer2012or} used a partially observable Markov decision process (POMDP) to take individualized mammography screening decisions while some personal risk features in addition to age and screening history of each patient is considered. They show that by considering these strategies the number of false positive results decreased and quality-adjusted-life-years (QALYs) are improved. This model is extended by \cite{cevik2018analysis} and resource constraints are added. They show that allocating capacity efficiently among individuals with different cancer risk levels leads to significant QALYs gains.

Another important application of MDPs in treatment of the diseases is liver and kidney transplant decisions \cite{alagoz2004optimal, alagoz2005incorporating, segev2005kidney, alagoz2007choosing}. \cite{alagoz2004optimal} created an infinite horizon MDP to determine when a patient with end-stage liver disease such as hepatitis C accepts a living-donor transplant. Depending on the quality of the match with the donor and the current health status of the patient, the model determines whether the transplant increases the expected total lifetime of the patient and whether the transplant should be done.

Our contributions to the literature are as follows: 1) We offer a mathematical model in the area of TB screening, a first in the relevant literature. All previous studies apply simulation or cost evaluation methods that compare the cost of different screening methods. In other words, there is no optimization model to determine TB screening plans. 2) To the best of our knowledge, our study is the first that categorize healthcare employees based on factors that affect the results of TB screening and finds the best test for each group to minimize the expected total cost by considering the infection rate of healthcare employees. 3) Our study is the first to develop a discrete time, infinite horizon MDP model to determine the optimal time between tests for each healthcare employee group in addition to detecting the best test type for each group. 4) For the first time in the literature, we consider the potential impacts and costs of a false negative result (i.e., probable spread of disease) in our formulation. 

\section{An MDP Model for the TB Test Scheduling Problem}
We propose a discrete time infinite horizon MDP model to formulate the problem. The decision epochs, the state space, the action set, transition probabilities, and cost parameters are described below. 

\subsection{Decision Epochs}
The decisions on whether or not an employee group should take a TB test and if so what type of test they should take are made annually. Each year, new employees who started their work in the hospital have to take either the blood or the skin test. Current employees with no TB history are also eligible to take a TB test.

\subsection{The State Space}
The employees are classified based on their salaries and risk of infection. Considering risk for classification purposes is perhaps obvious, but considering salary groups may not be. Unlike the other studies in the literature, our model captures the opportunity cost of lost time by employees in deciding which test to administer and how often. Thus, considering salary is important in capturing this opportunity cost. Let $\mathcal{I}=\{1,...,I\}$ be the set of employee types based on salary, and $\mathcal{J}=\{1,...,J\}$ be the set of employee types based on infection risk. The infection risk groups are categorized based on the employees' work locations, the specific job they do, and their BCG vaccination history. 

The state of the system is determined by the number of new and ongoing employees with no positive TB history, because the result of the test for employees with positive TB  history will almost always be positive. In Prisma Health, these employees with positive TB history fill up a questionnaire each year, and a decision regarding the necessity of an X-Ray test is made based on their responses. We let $y^{t}_{ij}$ be the number of current (not new) employees of salary group $i \in \mathcal{I}$ and risk group $j \in \mathcal{J}$ who are still employed at the healthcare facility during year $t$. The new employees in salary group $i$ and risk group $j$ who join the health systems in year $t$ are represented by $x^{t}_{ij}$.

Employees in different risk groups have different rates of TB infection. The rate of infection of each group in each year depends on the undetected infected employees in each of the groups in the previous year. In particular, undetected infected employees spread TB among other employees and increase the infection rate. Thus, we need to keep track of the number of undetected infected employees. Thus, let $u^t_{ij}$ be the number of undetected infected employees of salary group $i$ and risk group $j$ in the beginning of year $t$.
Thus, the state space takes the form
\begin{align}
{\vec{s}^t} = (\vec{x}, \vec{y}, \vec{u}) = (x^t_{ij}; y^t_{ij};u^{t}_{ij}), \ \ i\in \mathcal{I}; j\in \mathcal{J}
\end{align}
 We let $M^x_{ij}$, $M^y_{ij}$ and $M^u_{ij}$ be the maximum value for current number of employees, new arrivals, and infected employees of salary group $i$ and risk group $j$, respectively. Thus, The state space is finite since $x^t_{ij}$, $y^t_{ij}$, and $u^{t}_{ij}$ are bounded.

\subsection{The Action Set}
The actions that are taken in each state determine whether to administer a test and type of the test to be administered for each employee group. Recall that new employees have to take a test in their first year of employment. Let $a_{sij}^{xt}$ and $a_{bij}^{xt}$ denote the decisions regarding the type of the test for new employees of salary group $i$ and risk group $j$ in year $t$ where $a_{sij}^{xt}=1$ if skin test is selected and $0$ otherwise. Similarly, $a_{bij}^{xt}=1$ if blood test is selected and $0$ otherwise. Clearly, one of the tests must be selected for new employees. Thus,
\begin{align}
a_{sij}^{xt} + a_{bij}^{xt} =1, \label{SS1:conss1} \ \ i \in \mathcal{I}; j \in \mathcal{J}.
\end{align}  
Let $a_{sij}^{yt}$ and $a_{bij}^{yt}$ be the decisions regarding the type of the test for current employees of salary group $i$ and risk group $j$ at time $t$, where $a_{sij}^{yt} = 1 $ if skin test is selected and $0$ otherwise. Similarly, $a_{bij}^{yt}=1$ if blood test is selected and $0$ otherwise. Obviously, if both $a_{sij}^{yt}$ and $a_{bij}^{yt}$ are zero, the employees of salary group $i$ and risk group $j$ do not take the test at time $t$. Thus,
\begin{align}
0 \leq a_{sij}^{yt} + a_{bij}^{yt} \leq 1, \label{SS2:conss2} \ \ i \in \mathcal{I}; j \in \mathcal{J}.
\end{align}  
The action set takes the form
\begin{align}
\vec{a}_s^t = (a^{xt}_{sij},a^{xt}_{bij},a^{yt}_{sij},a^{yt}_{bij}), \ \ i \in \mathcal{I}; j \in \mathcal{J}, 
\end{align}     
and the set of allowable actions in each state $\vec{s} \in S$ that satisfy constraints (\ref{SS1:conss1}) and (\ref{SS2:conss2}) are denoted by $A_{\vec{s}}$.

\subsection{Transitions}
Since a portion of employees stop working at the healthcare system during each year, we define $l_{ij}^{t}$ as the number of employees of salary group $i$ and risk group $j$ who leave the healthcare facility in year $t$, which is defined as $l_{ij}^{t} \sim Binomial(y^{t}_{ij},p_{ij}^l)$, where $p_{ij}^l$ is the probability that an employee of salary group $i$ and risk group $j$ leaves the system. Thus, $l_{ij}^{t}$ represents the number of employees with no positive TB history who leave.

The new employees who are selected to take skin test, have to take the test twice. The second test must be taken within 7-21 days of the first test.
Let $d_{sij}^{xt}$ and $d_{sij}^{yt}$ be the total number of new and ongoing employees of salary group $i$ and risk group $j$ who take the skin test in year $t$. Since there is no difference between the testing protocols for new and ongoing employees if they are selected to take a blood test, we let $d_{bij}^{t}$ be the total number of employees of salary group $i$ and risk group $j$ who take the blood test in year $t$. We define 
\begin{align}
&d_{sij}^{xt} = a_{sij}^{xt}x^{t}_{ij}, \ \ i=1,...,I; j=1,...,J,\\
&d_{sij}^{yt} =  a_{sij}^{yt}(y^{t}_{ij}-l^t_{ij}), \ \ i=1,...,I; j=1,...,J,\\ 
&d_{bij}^{t} = a_{bij}^{xt}x^{t}_{ij} + a_{bij}^{yt}(y^{t}_{ij}-l^t_{ij})\ \ i=1,...,I; j=1,...,J. 
\end{align}
Each year, a random number new infections occur in each employee group. Let $\alpha_{ij}^t$ be the infection probability for an employee of salary group $i$ and risk group $j$ in year $t$. The infection probability of each group in a year depends on the proportion of undetected infected employees in all employee groups in the previous year, the likelihood that the employees of the group contact with employees of different groups, and the transmission probability conditional on such contact. The infection probability also depends on the percentage of TB infected patients who visit the healthcare facility, probability that an employee of this group contacts an a patient and the transmission probability conditional on such contact. Let $u_{ij}^t$ denote the total number of undetected infected employees of salary group $i$ and risk group $j$ in year $t$. We define $\rho_{ij,i'j'}$ as the probability that an employee of salary group $i$ and risk group $j$ contacts individuals of salary group $i'$ and risk group $j'$ during that year, and $\xi_{ij}$ as the TB transmission probability of an employee of salary group $i$ and risk group $j$ in case of contacting an infected individual. We define $\beta$ as the proportion of infected patients who visit the healthcare facility, $\nu_{ij}$ as the probability that an employee of salary group $i$ and risk group $j$ contacts a patient. Thus, we define $\alpha_{ij}^t$ as 
\begin{align}
\alpha_{ij}^{t} = \sum_{i'=1}^{I}{\sum_{j'=1}^{J}{  \rho_{ij,i'j'} \xi_{ij}\frac{u_{i'j'}^{t-1}}{x^{t-1}_{i'j'}+y^{t-1}_{i'j'}} }} +  \beta\nu_{ij} \xi_{ij} \ \ i=1,...,I; j=1,...,J
\end{align}
We define $n^{xt}_{sij}$ and $n^{yt}_{sij}$ as the number of infected new and ongoing employees of salary group $i$ and risk group $j$ who have to take the skin testin year $t$. We also define $n^t_{bij}$ as the number of infected employees of salary group $i$ and risk group $j$ who take the blood test at time $t$. We let $n^{xt}_{sij}$, $n^{yt}_{sij}$ and  $n^t_{bij}$ follow binomial distributions in the forms of $Binomial(d_{sij}^{xt},\alpha_{ij}^t)$, $Binomial(d_{sij}^{yt},\alpha_{ij}^t)$ and $ Binomial(d_{bij}^{t},\alpha_{ij}^t)$, respectively.

Due to possible false negative test results, a portion of infected employees are undetected. The probability of getting a false negative result depends on the BCG vaccination history of employees. Since BCG vaccination is one of the factors that we consider in defining risk groups, the probability of getting a false negative result depends on the employees' risk groups. Let $p_{sj}^{n}$ and $p_{bj}^{n}$ be the probability of getting a false negative result for an employee of risk group $j$ in skin test and blood test, respectively. The number of false negative results in skin test for ongoing employees of salary group $i$ and risk group $j$ at time $t$, $u_{sij}^{yt}$, follows binomial distribution in the form of $Binomial(n_{sij}^{yt},p_{sj}^{n})$. For new employees who get tested with a skin test, getting a false negative result means getting false negative results in both the initial and the follow up tests. Thus, the number of false negative results $u_{sij}^{xt}$, follows a distribution in the form of $Binomial(Binomial(n_{sij}^{xt},p_{sj}^{n}),p_{sj}^{n})$. We let $u_{bij}^{t}$ be the number of false negative results for employees of salary group $i$ and risk group $j$ in year $t$ who have taken blood test which follows binomial distribution $Binomial(n^t_{bij},p_{bj}^{n})$.

The number of new and ongoing employees of salary group $i$ and risk group $j$ who have got true positive skin test results are defined as $q^{xt}_{sij}$ and $q^{yt}_{sij}$, respectively, and calculated as $q^{xt}_{sij} = n_{sij}^{xt} - u_{sij}^{xt}$ and $q^{yt}_{sij} = n_{sij}^{yt} - u_{sij}^{yt}$. Similarly, $q^{t}_{bij}$ represents the total number of employees of salary group $i$ and risk group $j$ who have true positive blood test results in year $t$ where $q^{t}_{bij} = n_{bij}^{t} - u_{bij}^{t}$.  

According to our definitions, sum of $u_{sij}^{xt}$, $u_{sij}^{yt}$ and $u_{bij}^{t}$ shows the number of undetected infected employees in salary group $i$ and risk group $j$ that were given one of the TB tests in year $t$. There might also be infected employees in the employee groups that were not given any of the tests in year $t$. Let $d^t_{nij}$ be the number of employees in salary group $i$ and risk group $j$ that was not tested in year $t$, which is defined as    
\begin{align}
&d^t_{nij} = (1-a_{sij}^{yt}-a_{bij}^{yt})(y^{t}_{ij}-l^t_{ij}), \ \ i=1,...,I; j=1,...,J.
\end{align}
The number of infected employees in these groups of employees also follows binomial distribution with infection probability of each group. Let $u^t_{nij}$ be the number of infected employees of salary group $i$ and risk group $j$ that was not tested in year $t$ where, $u_{nij}^t \sim Binomial(d^t_{nij},\alpha_{ij})$. We define $u^t_{ij}$ as the total number of undetected infected employees where
\begin{align}
&u^t_{ij} = u_{sij}^{xt} + u_{sij}^{yt} + u_{bij}^t + u_{nij}^t, \ \ i=1,...,I; j=1,...,J.
\end{align}
Some uninfected employees might receive false positive results from the TB tests, which leads to a follow up X-ray at extra cost. We let $r^{xt}_{sij}$ and $r^{yt}_{sij}$ be the number of new and ongoing employees of salary group $i$ and risk group $j$ who have received false positive skin test results. Similarly, we let $r^{t}_{bij}$ show the number of employees of salary group $i$ and risk group $j$ who have got false positive blood test results. These variables follow the distributions $r^{xt}_{sij}\sim Binomial(d_{sij}^{xt} - n_{sij}^{xt},p^p_{sj})$ + $Binomial(d_{sij}^{xt} - n_{sij}^{xt} - Binomial(d_{sij}^{xt} - n_{sij}^{xt},p^p_{sj}),p^p_{sj})$, $r^{yt}_{sij}\sim Binomial(d_{sij}^{yt} - n_{sij}^{yt},p^p_{sj})$ and $r^t_{bij}\sim Binomial(d_{bij}^{t} - n_{bij}^{t},p^p_{bj})$ where $p^p_{sj}$ and $p^p_{bj}$ are the probabilities that an employee of salary group $i$ and risk group $j$ gets a false positive result in skin and blood tests, respectively.

Once the testing decisions are made, the stochastic elements that determine the state transition are the new arrivals of employees, the employees who left the system, and the total number of infections. The evolution of state space elements are captured by the following equations:
\begin{align}
&y^{t+1}_{ij} = y^{t}_{ij} + x^{t}_{ij} - l_{ij}^{t} - n_{sij}^{xt} - n_{sij}^{yt} - n_{bij}^{t} - u_{nij}^{t}, \ \ i=1,...,I; j=1,...,J\\
&u^{t+1}_{ij} = u_{sij}^{xt} + u_{sij}^{yt} + u_{bij}^t + u_{nij}^t, \ \ i=1,...,I; j=1,...,J\ \  
\end{align}
The trajectory of the system is represented by $\{ (\vec{s}^{t},\vec{a}_s^t): t=1,2,... \}$, where $\vec{s}^t$ is the state of the system and $\vec{a}_s^t$ is the action that is taken in year $t$. The stochastic evolution of the system is represented by $\vec{s}^{t+1}=F(\vec{s}^t,\vec{a}_s^t,g(\vec{s}^t,\vec{a}_s^t))$, where $F(.,.,.)$ is a transfer mapping and $g(\vec{s}^t,\vec{a}_s^t)$ is a random element that contains all the random quantities in the system at time $t$. These definitions are used for defining the value function.
\subsection{The costs}
The cost associated with each state-action pair drives from four sources: cost of doing the tests $(c^b, c^s)$, cost of doing the X-ray $(c^x)$, cost of lost time of employees $(c^l_i)$, and cost of undetected infections $(c^u_i)$.
\begin{align}
c(\vec{s},\vec{a}) = \sum_{i=1}^{I}{ \sum_{j=1}^{J}{\Big( c^b d_{bij}^{t} + c^s d_{sij}^{t} + c^x (q^{xt}_{sij}+q^{yt}_{sij}+q^t_{bij}+r^{xt}_{sij}+r^{yt}_{sij}+r^t_{bij}) + c^u_{ij}  u^{t}_{ij} + c^l_i w^t_{ij} } \Big)},
\end{align} 

where $w^t_{ij}$ is the expected total time spent at the testing clinic for employees of salary group $i$ and risk group $j$ in year $t$. For employees who are taking blood test in year $t$, $w^t_{ij}$ either only depends on the time spent administering the blood test, or also includes the time spent on X-Ray if result of the blood test is positive. However, for employees who are taking the skin test, the total skin test time depends on whether or not the employee should take the test twice or once (i.e., employee is new or not). There is also a possibility that the employee is not able to complete the second step of the skin test within the required time and has to repeat both steps. This also affects the expected total time an employee spends in taking the skin test.  

\subsection{Optimality equation}
The value function $v(\vec{s})$ corresponds to the total expected discounted cost for state $\vec{s}$ over the infinite horizon. 
\begin{align} 
&v(\vec{s}) = \min_{\vec{a} \in A_{\vec{s}}} \bigg\{ c(\vec{s},\vec{a})+\lambda \mathbb{E} \big(v(\vec{s'})\big) \bigg\}, \ \ \forall \vec{s} \in S, \label{eq:bellman}
\end{align}
where the expectation is taken with respect to $s'=F(\vec{s}^t,\vec{a},g(\vec{s}^t,\vec{a}))$ and $\lambda \in [0,1)$ is a discount factor. To find the optimal policy we need to solve equation (\ref{eq:bellman}). Since the state space and action set are finite, there exists a stationary optimal policy. However, the size of the state space make a direct solution to (\ref{eq:bellman}) impractical.  

\section{Approximate Dynamic Programming} 
Due to the curse of dimensionality, we use approximate dynamic programming to estimate the optimal policy. First, we transform the MDP model into its equivalent linear program (LP) as follows:
\begin{align}
 &\max \sum_{\vec{s} \in S}^{}{\gamma(\vec{s})v(\vec{s})} \\
 &\nonumber\text{s.t.} \\
 &c(\vec{s},\vec{a}) + \lambda \sum_{\vec{s'} \in S}^{}{\mathbb{P}(\vec{s'}|\vec{s},\vec{a})v(\vec{s'})} \geq v(\vec{s}) \ \ \forall \vec{s} \in S, \vec{a} \in A_{\vec{s}},
\end{align}
where $\mathbb{P}(.|.,.)$ is the transition probability and $\gamma(.)$ is the probability distribution over the initial state of the system. The LP formulation does not avoid the curse of dimensionality. Thus, we approximate the value function by using a specific parameterized form where the interactions of employees from different groups are not considered. The resulting approximate value function can be written as the summation of separate value functions for each employee group as follows:
\begin{align}
&v(\vec{s}) \approx \sum_{i=1}^{I}{\sum_{j=1}^{J}{v_{ij}(\vec{s}_{ij})}}, \ \ \forall \vec{s} \in S,
\end{align}
where $\vec{s}_{ij}$ represents the state of the system for employees of salary group $i$ and risk group $j$ which shows the number of current employees, new arrivals and infected employees of salary group $i$ and risk group $j$.  Moreover, $v_{ij}(\vec{s}_{ij})$ denotes the value function for employees of salary group $i$ and risk group $j$ that is defined as follows:
\begin{align}
&\nonumber v_{ij}(\vec{s}_{ij}) = \min_{\vec{a}_{ij} \in A_{\vec{s}_{ij}}} \bigg\{ c^b d_{bij}^{t} + c^s d_{sij}^{t} + c^x (q^{xt}_{sij}+q^{yt}_{sij}+q^t_{bij}+r^{xt}_{sij}+r^{yt}_{sij}+r^t_{bij}) + c^u_i  u^{t}_{ij} + c^l_i w^t_{ij} + \\
& \lambda \Big( \sum_{\vec{s'_{ij}} \in S_{ij}}^{}{\mathbb{P}(\vec{s'_{ij}}|\vec{a_{ij}},\vec{s_{ij}}) \times v_{ij}(\vec{s'}_{ij})} \Big) \bigg\}, \ \ \forall i=1,...,I, j=1,...,J, \vec{s}_{ij} \in S_{ij}.
\end{align}

When the LP-based ADP presented here is used, the number of decision variables grows linearly with the number of employees in each group. In contrast, if the original LP formulation is used, the number of variables grows exponentially with the number of employees in each group. Thus, we expect that this approximate reformulation will help addressing the computational challenges due to curse of dimensionality.

Since $w^t_{ij}$ includes both the service time (time of taking the tests) and waiting time of employees of salary group $i$ and risk group $j$ at time $t$, states of other groups affect the value of $w^t_{ij}$. In the LP-based ADP, we estimate $w^t_{ij}$ for each group independently of the states of the other employee groups. Thus, we calculate an upper bound for $w^t_{ij}$ and use this upper bound in our approximation.  To compute the upper bound, we assume that all others groups have the largest possible number of employees that we can have in each group, and all employees are taking the skin test. We build a pre-processing simulation model to estimate $w^t_{ij}$ for each group depending on its current state.

Recall that we have defined $\gamma(\vec{s})$ $\forall \vec{s} \in S$ as the probability distribution over the initial state of the system. We let $\gamma_{ij}(\vec{s}_{ij})$ $\forall \vec{s}_{ij} \in S_{ij}$ be the marginal initial state distribution for employees of salary group $i$ and risk group $j$, where $S_{ij}$ is the set of all possible states for each group. Also, let $\Omega =\{ (\vec{s}_{ij},\vec{a}_{ij}) : \vec{s}_{ij} \in S_{ij}, \vec{a}_{ij} \in A_{\vec{s}_{ij}}, \forall i=1,...,I, j=1,...,J \}$ be the set of all feasible state-action pairs.
%The LP-based ADP can be expressed as
%\begin{align}
%&\max \sum_{i=1}^{I}{\sum_{j=1}^{J}{\sum_{\vec{s}_{ij} \in S_{ij}}^{}{\gamma_{ij}(\vec{s}_{ij})v_{ij}(\vec{s}_{ij})}}} \\
%&\nonumber\text{s.t.} \\
%&c(\vec{s}_{ij},\vec{a}_{ij}) + \lambda \sum_{\vec{s'}_{ij} \in S_{ij}}^{}{\mathbb{P}(\vec{s'}_{ij}|\vec{s}_{ij},\vec{a}_{ij})v(\vec{s'}_{ij})} \geq v(\vec{s}_{ij}) \ \ \forall i=1,...,I, j=1,...,J, \vec{s}_{ij} \in S_{ij}, \vec{a}_{ij} \in A_{\vec{s}_{ij}} \\
%&v_{ij}(\vec{s}_{ij}) \geq 0, \ \ \forall i=1,...,I, j=1,...,J, \vec{s}_{ij} \in S_{ij}
%\end{align}
The dual of the LP-based ADP model can be expressed as
\begin{align}
&\min \sum_{i=1}^{I}{\sum_{j=1}^{J}{\sum_{(\vec{s}_{ij},\vec{a}_{ij}) \in \Omega}^{}{c(\vec{s}_{ij},\vec{a}_{ij})\delta(\vec{s}_{ij},\vec{a}_{ij})}}} \label{eq:c3}\\
&\nonumber\text{s.t.} \\
&\sum_{\substack{(\vec{s'}_{ij},\vec{a}_{ij}) \in \Omega \\ \vec{s'}_{ij}=\vec{s}_{ij}}}^{}{\delta(\vec{s'}_{ij},\vec{a}_{ij})}  - \lambda \sum_{(\vec{s'}_{ij},\vec{a}_{ij}) \in \Omega}^{}{\mathbb{P}(\vec{s}_{ij}|\vec{s'}_{ij},\vec{a}_{ij}) \delta(\vec{s'}_{ij},\vec{a}_{ij})} \geq \gamma_{ij}(\vec{s}_{ij}), \ \ \forall i=1,...,I, j=1,...,J, \vec{s}_{ij} \in S_{ij} \label{eq:c1} \\ 
&\delta(\vec{s}_{ij},\vec{a}_{ij}) \geq 0, \ \ \forall i=1,...,I, j=1,...,J, (\vec{s}_{ij},\vec{a}_{ij}) \in \Omega  \label{eq:c2}
\end{align}
In the above formulation, $c(\vec{s}_{ij},\vec{a}_{ij})$ shows the immediate cost for employees of salary group $i$ and risk group $j$ in state $\vec{s}_{ij}$ while action $\vec{a}_{ij}$ is taken.
The advantage of solving the dual problem is that the number of constraints is fewer. However, the number of variables is still very large. Thus, we use the column generation algorithm to obtain the optimal solution. The column generation algorithm is started with a small set of feasible state-action pairs (i.e., columns) to the dual problem, which is called the master problem. Then one or more violated constraints in the primal problem are found by solving a subproblem. The state-action pair(s) corresponding to these violated constraints are added to the master problem as new columns. The procedure continues until no primal constraint is violated.

We consider a master problem associated with the formulation (\ref{eq:c3}-\ref{eq:c2}), and let $v^D_{ij}(\vec{s}_{ij})$ be the dual variable associated with constraint (\ref{eq:c1}) for all $\vec{s}_{ij} \in S_{ij}$ for $i=1,...,I, j=1,...,J$. The corresponding subproblem can be written as:
\begin{align}
&z^{sub} = \min_{(\vec{s}_{ij},\vec{a}_{ij}) \in \Omega} \bigg\{ \sum_{i=1}^{I}{\sum_{j=1}^{J}{\big(c(\vec{s}_{ij},\vec{a}_{ij}) - v^D_{ij}(\vec{s}_{ij}) + \lambda \sum_{\vec{s'}_{ij} \in S_{ij}}^{} {\mathbb{P}_{ij}(\vec{s'}_{ij}|\vec{s}_{ij},\vec{a}_{ij}) v^D_{ij}(\vec{s'}_{ij}) }\big)}}. \bigg\} \label{S1:sub}
\end{align}
 Formulation (\ref{S1:sub}) is a generalized multidimensional knapsack problem. To solve the subproblem, we reformulate it as an integer program that is easy to implement in optimization solvers.

For each employee of salary group $i$ and risk group $j$ we let $\Omega_{ij} =\{ (\vec{s}_{ij},\vec{a}_{ij}) : \vec{s}_{ij} \in S_{ij}, \vec{a}_{ij} \in A_{\vec{s}_{ij}} \}$ be set of all feasible state-action pairs in the master problem and note that $f_{ij}(\vec{s}_{ij},\vec{a}_{ij})=c(\vec{s}_{ij},\vec{a}_{ij}) - v^D_{ij}(\vec{s}_{ij}) + \lambda \sum_{\vec{s'}_{ij} \in S_{ij}}^{} {\mathbb{P}_{ij}(\vec{s'}_{ij}|\vec{s}_{ij},\vec{a}_{ij}) v^D_{ij}(\vec{s'}_{ij})}$ is the corresponding reduced cost for such a state-action pair. Thus, the subproblem finds the feasible state-action pairs such that $\sum_{i=1}^{I}{\sum_{j=1}^{J}{f_{ij}(\vec{s}_{ij},\vec{a}_{ij})}} < 0$. We let  $z_{ij}^{\vec{s}_{ij},\vec{a}_{ij}}$ be a binary variable that is one if action $\vec{a}_{ij}$ is chosen in state $\vec{s}_{ij}$ for employees of salary group $i$ and risk group $j$; and zero, otherwise. Then, we transform the subproblem into a generalized assignment problem which chooses one feasible state -action pair for each employee group. The formulation for the assignment problem is as follows:
\begin{align}
&\min \sum_{i=1}^{I}{\sum_{j=1}^{J}{\sum_{(\vec{s}_{ij},\vec{a}_{ij}) \in \Omega_{ij}}^{}{f_{ij}(\vec{s}_{ij},\vec{a}_{ij})z_{ij}^{\vec{s}_{ij},\vec{a}_{ij}}}}} \label{eq:assignment} \\ 
&\nonumber\text{subject} \ \ \text{to} \\
&\sum_{(\vec{s}_{ij},\vec{a}_{ij}) \in \Omega_{ij}}^{}{z_{ij}^{\vec{s}_{ij},\vec{a}_{ij}} =1}, \ \ i=1,...,I, j=1,...,J \\
&z_{ij}^{\vec{s}_{ij},\vec{a}_{ij}} \in \{ 0,1\}, \ \  i=1,...,I, j=1,...,J \ \ \forall (\vec{s}_{ij},\vec{a}_{ij}) \in \Omega_{ij}
\end{align}
The objective function (\ref{eq:assignment}) chooses the column with minimum reduced cost. This new formulation can be implemented and solved using commercial optimization solvers.  

\section{Numerical Study and Results}
In this section, we implement the approximate dynamic programming scheme described in the previous section using data obtained from one of the Prisma Health hospitals and we obtain the optimal TB screening policy. Based on the optimal policy we propose some benchmark policies and simulate them to show how these policies improve the system compare to the current policies.
\subsection{Input data}
In this study we consider three employee salary groups; (i) physicians, (ii) nurses and (iii) other employees, and three employee risk groups; (i) BCG vaccinated employees, (ii) employees who work in high risk locations, (iii) and employees who work in low risk locations. Thus, a total of nine groups are presumed. We need to mention that the yearly arrivals of employees are assumed to be stochastic and modeled as truncated Poisson processes. Yearly arrival rates, leaving probabilities, probabilities of contacting patients and employees, TB transmission probabilities, TB infection probability, false positive and false negative probabilities, and cost parameters are the parameters that are used in the MDP model. The estimation for these parameters are presented Table \ref{table:parametersestimation}. 

Probabilities of contacting patients for different employee groups $\nu_{ij}$, are estimated based on our discussions with Prisma Health physicians. In using ADP, we assume that employees of different groups do not contact with each other and employees who are in the same group definitely contact each other. Thus, $\rho_{ij,i'j'}$, which is defined as the probability that an employee of salary group $i$ and risk group $j$ contacts individuals of salary group $i'$ and risk group $j'$, is zero if $i \neq i'$ or $j \neq j'$ and it is one for employees in the same group. We base our estimates of TB Transmission probabilities, $\xi_{ij}$, for each group, we on the study by \cite{ahmed2011molecular}. Since Prisma health did not have exact information about the percentage of patients with TB who come to the hospital, $\beta$, we used the information from \cite{mullie2017revisiting} to estimate this parameter. In that study, the average percentage of individuals infected with TB among healthcare workers is reported to be about $2 \%$. Thus, we simulated and calibrated the system under the hospital's current policy (i.e., annual skin tests for all employeees) so that the value of $\alpha_{ij}^t$ in the long run is $2 \%$. The results of this simulation showed that value of $\beta$ which brings $\% 2$ for $\alpha_{ij}^t$ is $0.1$. Hence, we use this value in our experiments. The false positive and false negative probabilities of the tests are taken from \cite{duchin1997comparison, king2015t, dorman2014interferon, jensen2005guidelines}. The remaining parameters used in the model (i.e., the cost parameters) are estimated based on the hospital data. 
\subsection{Experimental setup}
The optimization problem is implemented in C++. The subproblems and master problems are solved on an Intel Core i7-9700  CPU  utilizing  the  Gurobi  9.0  solver. The computational time of solving the total problem for one instance is about 120 minutes.

\subsection{Optimal policy for the base model}
\begin{figure}[ht!]
\minipage{0.45\textwidth}
  \includegraphics[width=\linewidth]{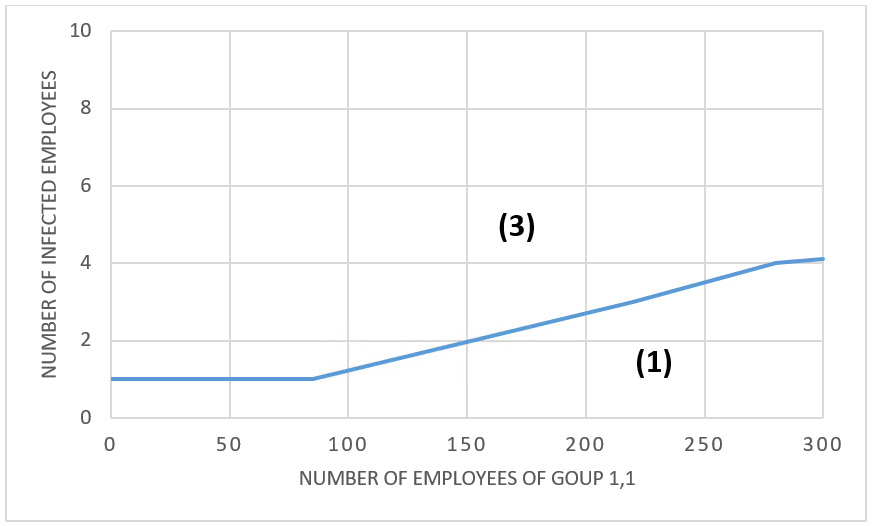}
  \subcaption{}
  \label{fig:group1a}
\endminipage\hfill
\minipage{0.45\textwidth}
  \includegraphics[width=\linewidth]{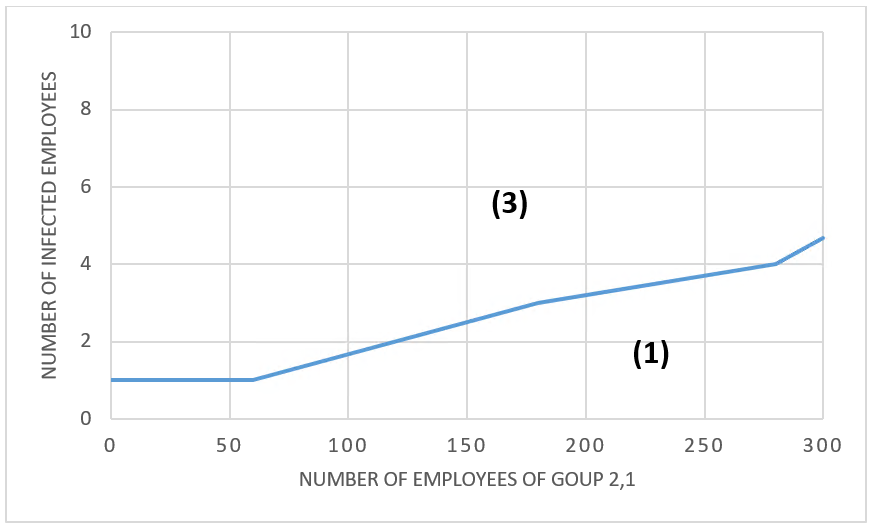}
  \subcaption{}
  \label{fig:group1b}
\endminipage\hfill
\begin{center}
\minipage{0.45\textwidth}%
  \includegraphics[width=\linewidth]{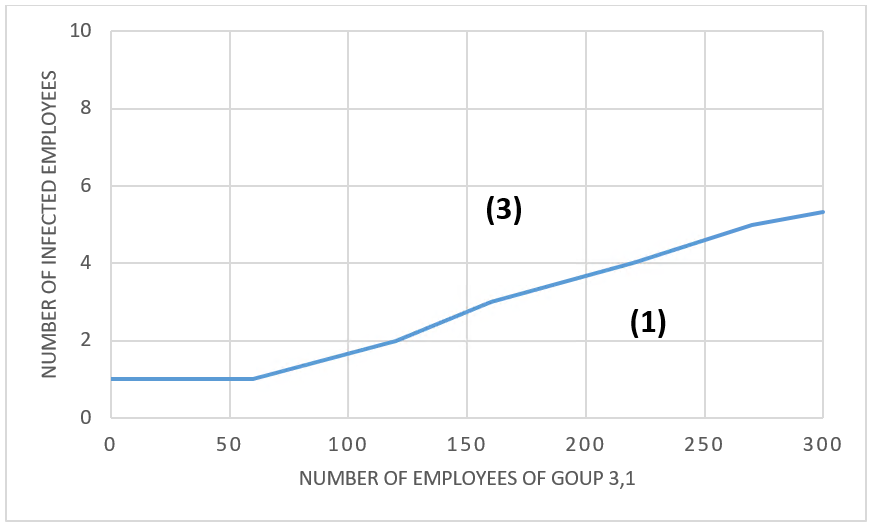}
  \subcaption{}
  \label{fig:group1c}
\endminipage
\end{center}
 \caption{Optimal policy for risk group 1}
\label{fig:group1}

\end{figure}

\begin{figure}[ht!]
\minipage{0.45\textwidth}
  \includegraphics[width=\linewidth]{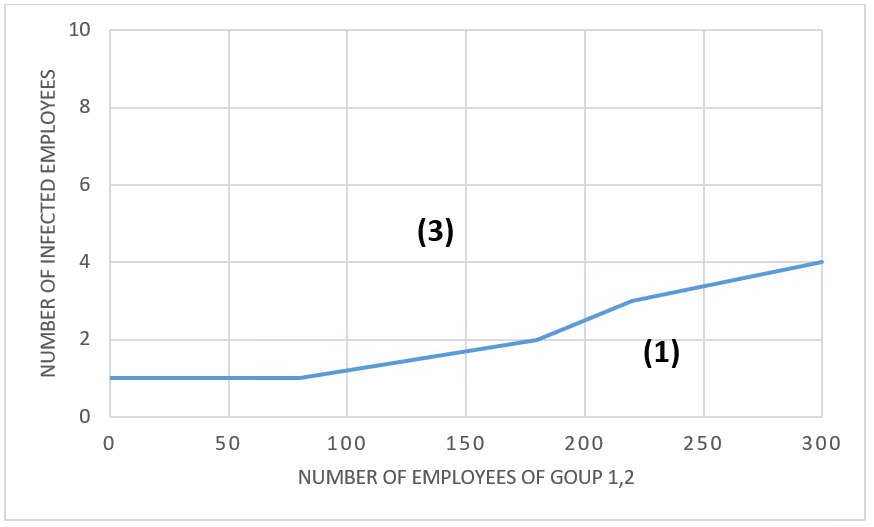}
  \subcaption{}
   \label{fig:group2a}
\endminipage\hfill
\minipage{0.45\textwidth}
  \includegraphics[width=\linewidth]{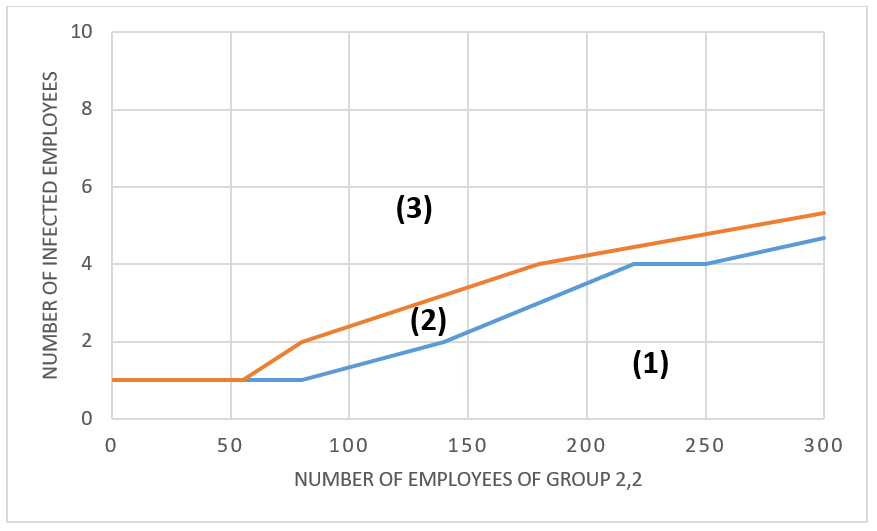}
  \subcaption{}
   \label{fig:group2b}
\endminipage\hfill
\begin{center}
\minipage{0.45\textwidth}%
  \includegraphics[width=\linewidth]{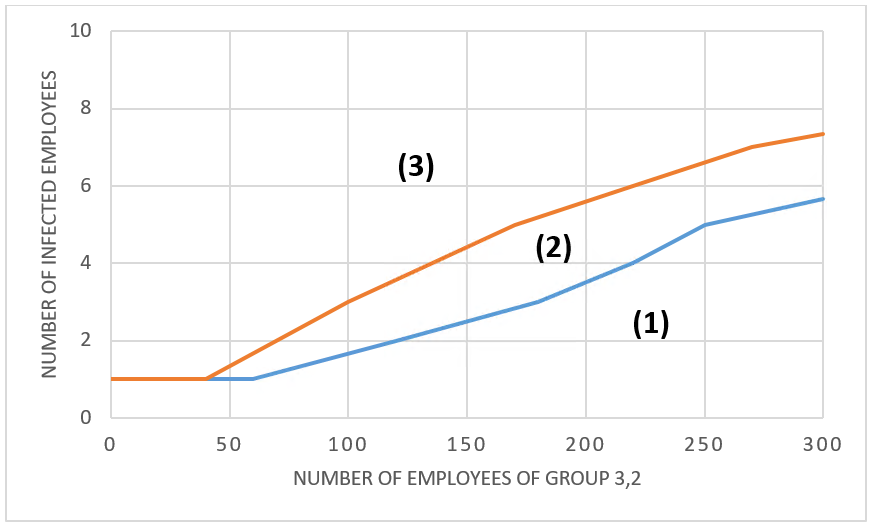}
  \subcaption{}
   \label{fig:group2c}
\endminipage
\end{center}
 \caption{Optimal policy for risk group 2}
\label{fig:group2}
\end{figure}

\begin{figure}[ht!]
\minipage{0.45\textwidth}
  \includegraphics[width=\linewidth]{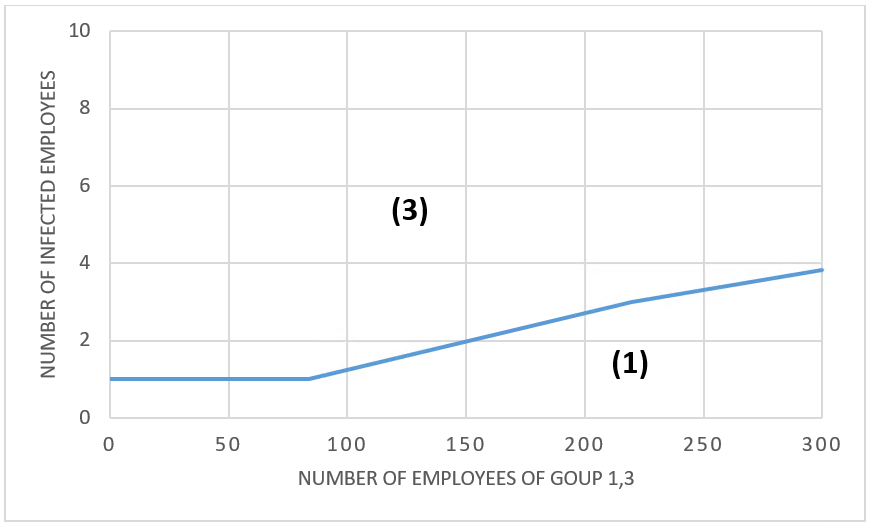}
  \subcaption{}
   \label{fig:group3a}
\endminipage\hfill
\minipage{0.45\textwidth}
  \includegraphics[width=\linewidth]{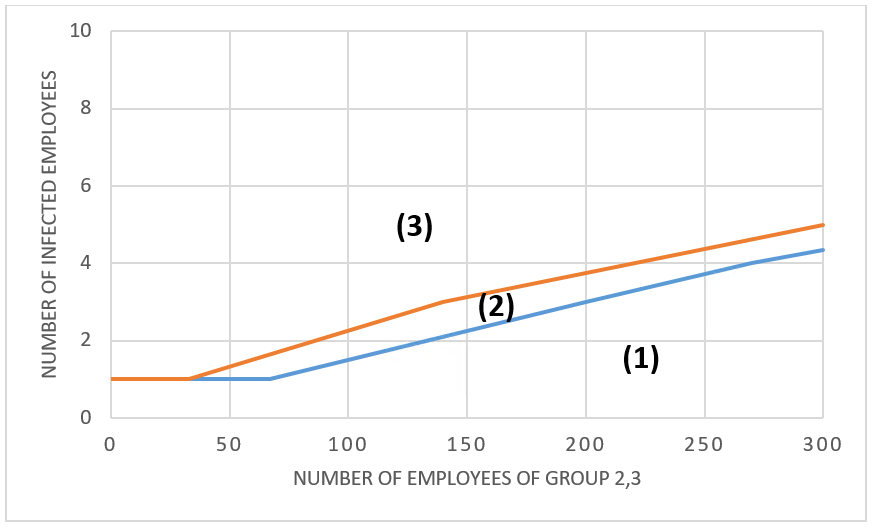}
  \subcaption{}
   \label{fig:group3b}
\endminipage\hfill
\begin{center}
\minipage{0.45\textwidth}%
  \includegraphics[width=\linewidth]{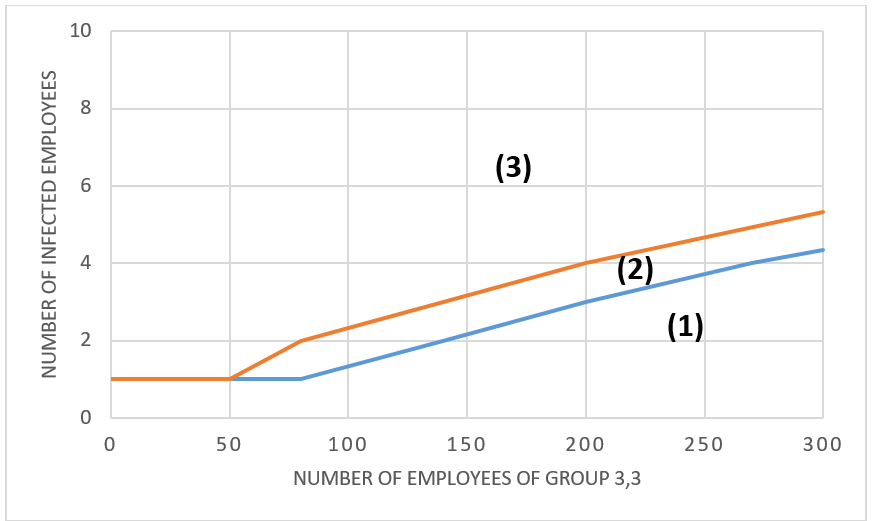}
  \subcaption{}
   \label{fig:group3c}
\endminipage
\end{center}
 \caption{Optimal policy for risk group 3}
\label{fig:group3}

\end{figure}

The optimal policy which dictates testing actions in each state is obtained by solving the ADP model. After analyzing the optimal policy, we observed that changes in number of new arrivals does not have huge impact on the optimal testing decision. Thus, to visualize the optimal policy and make its discussion easier, we created plots of current number of employees ($y^t_{ij}$) versus number of infected employees ($u^t_{ij}$). These plots are presented as Figures \ref{fig:group1}, \ref{fig:group2}, and \ref{fig:group3}. 

Each plot has regions that are numbered (1), (2) or (3). Region (1) shows the states where the optimal action is not taking any test. Regions (2) and (3) represent the states where the optimal action is taking skin test and blood test, respectively. In the following we briefly describe these graphs.

Figure \ref{fig:group1} shows the optimal policies for employees of risk group 1 (i.e., employees who had the BCG vaccination).  We observe that for employees across all  salary groups who had the BCG vaccination, the optimal decision is either no testing or the blood test depending on the level of infection spread among the employees. This is because the skin tests are less cost-effective for these employees due to the major false positive risk.

Figure \ref{fig:group2} depicts the optimal policies for employees of risk group 2 (i.e., employees who work at low risk locations). The results show that for physicians, the optimal decision is either no test or blood test. That is due to the high cost of lost time incurred when these employees get skin tests. On the other hand, for nurses and other employees depending on the number of employees and number of infected ones the decision would be either no test, skin test or blood test.

Finally, the optimal policies for risk group 3 (i.e.,  employees who work at high risk locations) is shown in Figure \ref{fig:group3}. As we expected, for physicians the optimal decision is either no test or blood test, and for nurses and others the optimal decision can be no test, skin test or blood test depending on the level of infection spread. 

\subsection{Simulation of the optimal policy}
The Optimal policy indicates the optimal action (no test, skin test, or blood test) for each employee group in different states. However, how often the tests are performed under the optimal policy is not readily known without further analysis. we are interested in knowing these frequencies since the current practice as well as CDC guidelines are based on testing frequency and not on the number of employees and infections detected. Thus, we simulate the optimal policy to estimate the time between tests for each employee group using the number of times that each group visits the states that require testing. The simulation model is implemented in C++ and ran for 100 years. The results of the simulation are presented in Table \ref{table:simulation}. We report the optimal action which is the type of the test administered and frequency of that type of test for each employee group. As mentioned above, the frequency of the tests is estimated based on the number of times that the test is administrated in a 100 year horizon..
\begin{table}[h]
\footnotesize
\centering
 \caption{Optimal actions based on the simulation}
  \label{table:simulation}
\begin{tabular}{|c|c|c|}
\hline
 Employee group & Optimal action for new employees  & Optimal action for current employees   \\ \hline \hline
1,1 & take blood test & take blood test every year \\ \hline
1,2 & take blood test & take blood test every year \\ \hline
1,3 & take blood test &  take blood test every year  \\  \hline
2,1 & take blood test & infrequently - take blood test if percentage of infected  \\ & & employees is greater than $1.6 \%$ \\ \hline
2,2 & take blood test & take blood test every 3 years \\  \hline
2,3 & take blood test & take blood test every 3 years\\  \hline
3,1 & take blood test & infrequently - take blood test if percentage of infected \\ & & employees is greater than $1.7 \%$ \\  \hline
3,2 & take blood test & infrequently - take blood test if percentage of infected\\ & & employees is greater than $2.2 \%$ and skin test if it \\ & & is between 1.7 and 2.2$\%$ \\  \hline
3,3 & take blood test & take blood test every 2 years \\ 
\hline
\end{tabular}
\end{table}

For some employee groups, the frequency of testing was very low, since the groups visited states that require testing infrequently. Thus, instead of just providing a recommendation of practically never testing these employees, we define the critical ratio between number of infected employees and current number of employees in these groups that would trigger testing. In other words, we define a threshold for percentage of infected employees where for larger percentages either a blood test or a skin test is required. These thresholds are extracted from the optimal policy.

In order to observe the advantages of the proposed policy given in Table \ref{table:simulation}, we simulate this policy and compare the result with the hospital's current policy which is administering skin tests every year for all employee groups. These simulations are also implemented in C++. In these simulations, we still assume that employees of different groups do not interact. Table \ref{table:simulation2} shows the corresponding results.
\begin{table}[h]
\footnotesize
\centering
 \caption{Comparison of the proposed policies and the current policy}
  \label{table:simulation2}
\begin{tabular}{|c|c|c|c|c|}
\hline
  & Current policy  & Proposed policy & Optimal policy \\ \hline \hline
Average yearly cost   & 17417 $\$$  &  8752 $\$$  & 6750 $\$$ \\ \hline
  Average infection rate   & 1.5 $\%$  &  1.5 $\%$ & 1.2 $\%$ \\
\hline
\end{tabular}
\end{table}
The proposed policy decreases the average screening cost of the healthcare facility by optimizing the type and frequency of the tests for different employee groups while the average infection rate does not increase. Based on theses results, it is beneficial for the hospital to update their screening guidelines for different employee groups instead of annually testing all employees regardless of risk and salary characteristics.
\section{Conclusion} \label{section:conclusion}
This study includes a dynamic TB screening model for healthcare employees. Although, rate of TB infection in the US is low, healthcare employees are still in risk of getting infected due to being in direct contact with patients. Thus, regular screening seems to be an effective preventive approach. For screening purposes, hospitals can administer skin or blood tests. These tests bring costs for the healthcare facility. Thus, providing an optimal screening policy that minimizes the infection rate and cost of the healthcare facility is important. 

This challenge motivated us to propose a dynamic model that specifies the best TB test and its frequency for each group of employee. We formulated this problem as a MDP and used ADP to estimate the solution. We used our formulation and a data set obtained from one of the Prisma Heealth hospitals to determine the optimal TB screening test for each employee group. Currently, an annual skin test is required in this hospital. Furtermore, we simulated the optimal policy to estimate the optimal screening frequency for each group of employees. Then, we simulated these resulting policy and compare the results with the current screening policy implemented at our partnering hospital. Results confirm the improvement in the average cost of the healthcare facility without an increase in the observed infection rate.

This study can be extended in multiple directions. Our model does not differentiate between the latent and active TB infections. However, latent and active TB may have different impacts on the results of the tests or may affect the spread of the disease differently. Furthermore, by using ADP we assumed that employees of different groups do not contact with each other. However, to consider a more realistic case, contacting with different groups can be considered.
\newpage
\section{Appendix}
\begin{table}[h]
\footnotesize
\centering
 \caption{Parameters values}
  \label{table:parametersestimation}
\begin{tabular}{|c|c||c|c||c|c|}
\hline
 Parameter& Value & Parameter & Value & Parameter & Value  \\ \hline
$\lambda_{11}$ & 4  & $\nu_{23}$ & 1 & $p^p_{s1}$ & 0.6 \\ 
$\lambda_{12}$ & 14 & $\nu_{31}$ & 0.75 & $p^p_{s2}$ & 0.27 \\
$\lambda_{13}$ & 10  & $\nu_{32}$ & 0.5 & $p^p_{s3}$ & 0.27 \\ 
$\lambda_{21}$ & 15 & $\nu_{33}$ & 1 & $p^n_{s1}$ & 0.04  \\ 
$\lambda_{22}$ & 50  & $\rho_{11,11}$ & 1 & $p^n_{s2}$ & 0.04 \\ 
$\lambda_{23}$ & 35  & $\rho_{12,12}$ & 1 & $p^n_{s3}$ & 0.04 \\ 
$\lambda_{31}$ & 4 & $\rho_{13,13}$ & 1 & $p^p_{b1}$ & 0.176 \\ 
$\lambda_{32}$ & 14  & $\rho_{21,21}$ & 1 & $p^p_{b2}$ & 0.176 \\ 
$\lambda_{33}$ & 10 &  $\rho_{22,22}$ & 1 & $p^p_{b3}$ & 0.176 \\  
$p^l_{11}$ & 0.15 & $\rho_{23,23}$ & 1 & $p^n_{b1}$ & 0.008  \\ 
$p^l_{12}$ & 0.15 & $\rho_{31,31}$ & 1 & $p^n_{b2}$ & 0.008 \\ 
$p^l_{13}$ & 0.15 & $\rho_{32,32}$ & 1 & $p^n_{b3}$ & 0.008\\ 
$p^l_{21}$ & 0.15 & $\rho_{33,33}$ & 1 & $c^b$ & 45 $\$$ \\ 
$p^l_{22}$ & 0.15 & $\xi_{11}$ & 0.05 &  $c^s$ & 8 $\$$\\ 
$p^l_{23}$ & 0.15 & $\xi_{12}$ & 0.22 & $c^x$ & 100 $\$$\\ 
$p^l_{31}$ & 0.15 & $\xi_{13}$ & 0.22 &  $c^l_1$ & 150 $\$/h$\\ 
$p^l_{32}$ & 0.15 & $\xi_{21}$ & 0.05  & $c^l_2$ & 30 $\$/h$\\ 
$p^l_{33}$ & 0.15 & $\xi_{22}$ & 0.22 & $c^l_3$ & 29 $\$/h$\\ 
$\nu_{11}$ & 1 & $\xi_{23}$ & 0.22 & $c^u_1$ & 5000 $\$$\\ 
$\nu_{12}$ & 1 &  $\xi_{31}$ & 0.05 & $c^u_2$ & 1000 $\$$\\ 
$\nu_{13}$ & 1 & $\xi_{32}$ & 0.22 & $c^u_3$ & 1000 $\$$ \\ 
$\nu_{21}$ & 1 & $\xi_{33}$ & 0.22 &  & \\ 
$\nu_{22}$ & 1 &$\beta$ & 0.1 &  & \\ 
\hline
\end{tabular}
\end{table}

\newpage
\bibliographystyle{plain}
\bibliography{main}

%%%%%%%%%%%%%%%%%
\end{document}